\documentstyle[12pt]{article}
\raggedright

\input amssym.def
\input amssym.tex

\newcommand{\eop}{\bigstar}  

\newcommand{\complB}{<\!\!\!\circ}
\newcommand{\card}[1]{{\vert #1 \vert} }
\newcommand{\otp}[1]{\hbox{otp($#1$)}}
\newcommand{\forces}{\Vdash}

\newcommand{\Dom}{{\rm Dom}}

\newcommand{\implies}{\Longrightarrow}

\newenvironment{Proof}{\noindent{\bf Proof.}}{\par\bigskip} 

\newtheorem{THEOREM}{Theorem}[section]

\newtheorem{Conclusion}[THEOREM]{Conclusion}

\newtheorem{LEMMA}[THEOREM]{Lemma}

\newtheorem{Main Theorem}[THEOREM]{Main Theorem}
\newenvironment{main Theorem}{\begin{Main Theorem}} 
{\end{Main Theorem}}

\newtheorem{Theorem}[THEOREM]{Theorem}

\newtheorem{Definition}[THEOREM]{Definition}

\newtheorem{Conventions}[THEOREM]{Conventions}

\newtheorem{Main Definition}[THEOREM]{Main Definition}
\newenvironment{main definition}{\begin{Main Definition}}
{\end{Main Definition}}

\newtheorem{Lemma}[THEOREM]{Lemma}

\newtheorem{Notation}[THEOREM]{Notation}

\newtheorem{Note}[THEOREM]{Note}

\newtheorem{Observation}[THEOREM]{Observation}

\newtheorem{Remark}[THEOREM]{Remark}

\newtheorem{Main Fact}[THEOREM]{Main Fact}
\newenvironment{main Fact}{\begin{Main Fact}}{\end{Main Fact}}

\newtheorem{Fact}[THEOREM]{Fact}

\newtheorem{Subfact}[THEOREM]{Subfact}

\newtheorem{Claim}[THEOREM]{Claim}

\newtheorem{Main Claim}[THEOREM]{Main Claim}
\newenvironment{main claim}{\begin{Main Claim}}{\end{Main Claim}}

\newtheorem{Corrolary}[THEOREM]{Corrolary}

\newtheorem{Subclaim}[THEOREM]{Subclaim}

\newtheorem{Corollary}[THEOREM]{Corollary}

\newtheorem{Proposition}[THEOREM]{Proposition}

\newtheorem{Discussion}[THEOREM]{Discussion}

\newenvironment{Proof of the Subfact}
{\noindent{\bf Proof of the Subfact.}}{\par\bigskip}

\newenvironment{Proof of the Fact}
{\noindent{\bf Proof of the Fact.}}{\par\bigskip}

\newenvironment{Proof of the Lemma}
{\noindent{\bf Proof of the Lemma.}}{\par\bigskip}

\newenvironment{Proof of the Claim}
{\noindent{\bf Proof of the Claim.}}{\par\bigskip}

\newenvironment{Proof of the Subclaim}
{\noindent{\bf Proof of the Subclaim.}}{\par\medskip}

\newenvironment{Proof of the Main Claim}
{\noindent{\bf Proof of the Main Claim.}}{\par\bigskip}

\catcode`\@=11
\def\@begintheorem#1#2{\rm \trivlist \item[\hskip \labelsep{\bf #1\
#2.}]}
\def\@opargbegintheorem#1#2#3{\rm \trivlist
      \item[\hskip \labelsep{\bf #1\ #2\ (#3).}]}
\catcode`\@=12

\newcommand{\stick}
{\mbox{{\hspace{0.4ex}$\hbox{\raisebox{1ex}{$\bullet$}}
\hspace*{-0.94ex}|$\hspace{0.6ex}}}}

\newcommand{\pr}{{\rm pr}}
\newcommand{\apr}{{\rm apr}}

\newcommand{\Max}{{\rm Max}}

\newcommand{\club}{\clubsuit}
\newcommand{\into}{\rightarrow}

\newcommand{\rest}{\upharpoonright}  

\newcommand{\acc}{\mathop{\rm acc}}
\newcommand{\satisfies}{\vDash}
\newcommand{\deq}{\buildrel{\rm def}\over =}

\newcommand{\HH}{{\cal H}}
\newcommand{\II}{{\cal I}}

\newcount\skewfactor
\def\mathunderaccent#1#2 {\let\theaccent#1\skewfactor#2
\mathpalette\putaccentunder}
\def\putaccentunder#1#2{\oalign{$#1#2$\crcr\hidewidth
\vbox to.2ex{\hbox{$#1\skew\skewfactor\theaccent{}$}\vss}\hidewidth}}
\def\name{\mathunderaccent\tilde-3 }

\title{Similar but not the same: various versions of $\clubsuit$ do not
coincide}
\author{Mirna D\v zamonja\\
Mathematics Department\\
University of Wisconsin--Madison \\
Madison, WI 53706, USA
\\
\scriptsize{dzamonja@math.wisc.edu}\and
Saharon Shelah\\
Mathematics Department\\
Hebrew University of Jerusalem\\
91904 Givat Ram, Israel\\
\scriptsize{shelah@sunset.huji.ac.il}
}

\date{}
\begin{document}

\baselineskip=16pt
\binoppenalty=10000
\relpenalty=10000
\raggedbottom

\maketitle

\begin{abstract}
We consider various versions of the $\clubsuit$ principle. This
principle is a
known consequence of $\diamondsuit$.
It is well known that $\diamondsuit$ is not sensitive to minor
changes in its definition, e.g. changing the guessing requirement
form ``guessing exactly" to ``guessing modulo a finite set".
We show
however, that this is not true for $\clubsuit$.
We consider some other variants of $\clubsuit$ as well.

\footnote{
This publication is numbered [F128]=[DjSh 574] in the list of
publications
of Saharon Shelah. 
The authors wish to thank the Basic Research Foundation
of the Israel Academy of Science for the support via their grant number
0327398.
In addition, Mirna D\v zamonja would like to thank the Hebrew
University of
Jerusalem and the Lady Davis Foundation for the
Forchheimer Postdoctoral Fellowship for the year 1994/95 when most of
the
research for this paper was done. The paper was distributed
in November 1995.
}
\end{abstract} 

\section{Introduction} In this paper we consider various natural
variants of
$\clubsuit$ principle. We answer questions of S. Fuchino and
M. Rajagopalan.

The principle was introduced by A. Ostaszewski
in
\cite{Ost}. It is easy to see that $\clubsuit$ follows
from $\diamondsuit$, and in fact it is true that $\diamondsuit$ is
equivalent to $\clubsuit+CH$, by an argument of K. Devlin presented in
\cite{Ost}. By ([Sh 98,\S5]) $\diamondsuit$ and $\clubsuit$ are not
equivalent,
that is, it is consistent to have $\clubsuit$ without having $CH$.
Subsequently J. Baumgartner,
in an unpublished note, gave an alternative proof, via a forcing which
does not
collapse $\aleph_1$ (unlike the forcing in [Sh 98]). P. Komj\'ath [Ko],
continuing the proof in
[Sh 98, \S 5] proved it consistent to have $MA$
for countable partial orderings $+\neg CH$, and $\club$. 
Then S. Fuchino, S. Shelah and  L. Soukup [FShS 544] proved
the same, without collapsing $\aleph_1$.

\indent The original R. Jensen's formulation of $\diamondsuit$
(\cite{Jensen})
is about the existence of a sequence $\langle
A_\delta:\,\delta<\omega_1\rangle$
such that every $A_\delta$ is an unbounded subset of $\delta$,
and for every $A\in [\omega_1]^{\aleph_1}$, we have
$A\cap\delta=A_\delta$
stationarily often. Many equivalent reformulations
can be obtained by using coding techniques (see \cite{Kunen}).
As a well known example, we mention K. Kunen's proof (\cite{Kunen})
that 
$\diamondsuit^-$ is equivalent to $\diamondsuit$. Here $\diamondsuit^-$
is 
the version of $\diamondsuit$ which says that there is 
a sequence $\langle
\{A^\delta_n:\,n<\omega\}:\,\delta<\omega_1\rangle$,
each $A_\delta^n\subseteq\delta$, and
for every $A\in [\omega_1]^{\aleph_1}$, we stationarily often have that
$A\cap\delta=A^\delta_n$ for some $n$. 

We consider
the question asking if $\clubsuit$ has a similar invariance property.
To be precise, we shall below formulate some versions of $\clubsuit$,
and ask
if any two of them are equivalent. We are particularly
interested in those versions of $\clubsuit$ which have the property that
the parallel version of $\diamondsuit$ is equivalent to
$\diamondsuit$.
The main result of the paper is that almost all of the
$\clubsuit$-equivalences we
considered,
are consistently false. 

Versions of $\clubsuit$ which are weaker than the
ones we consider,
are already known to be weaker
than $\clubsuit$. Namely, in his paper \cite{Ju}, I. Juh\'asz considers 
the principle $\clubsuit'$ claiming the existence of a
sequence $\langle\langle A^\delta_n:\,n<\omega\rangle:\,\delta\mbox{ limit }
<\omega_1\rangle$ where
for any $\delta$ sets $\{A^\delta_n:\,n<\omega\}$
are disjoint, and such that for every $A\in [\omega_1]^{\aleph_1}$ there
is $\delta$ such that for all $n$ we have $\sup(A^\delta_n\cap\omega_1)
=\delta$. I. Juh\'asz shows that $\clubsuit'$ is true in any extension
by a Cohen real. 

We heard of the question on the equivalence between $\clubsuit$
and $\clubsuit^\bullet$ from F. Tall, who heard it from
J. Baumgartner. J. Baumgartner credited the
question to F. Galvin, who credited it to M. Rajagopalan.
And indeed, M. Rajagopalan asked this question in \cite{Raja},
where he introduced $\clubsuit^\bullet$ (denoted there by
$\clubsuit_F$). In the same paper
M. Rajagopalan also introduced $\clubsuit^2$ (denoted there by
$\club^\infty$)
and showed that $CH+\clubsuit^2$
suffices for the Ostaszewski space.
He also asked if $\clubsuit^2$ was equivalent to $\clubsuit$.
The answer is negative by Theorem \ref{1 not 0} below.

Most of the other equivalence questions we consider here were first
asked by S.
Fuchino. 

We now proceed to give the relevant definitions.

\begin{Definition} We define the meaning of 
the principle $\clubsuit^l_\Upsilon$ for $l$
ranging in $\{0,1,2,\bullet\}$ and $\Upsilon$ a limit ordinal
$<\omega_1$.
(If $\Upsilon=\omega$ then we omit it from the notation.)
\smallskip

\underline{Case 1.} $l=0$

For some stationary set $S\subseteq\omega_1\cap LIM$,
there is a sequence $\langle A_\delta:\,\delta\in S
\rangle$ such that
\begin{description}
\item{{\em{(a)}}}
 $A_\delta$ is an unbounded subset of $\delta$.
\item {{\em (b)}} $\otp{A_\delta}=\Upsilon$.
\item{{\em (c)}} For every unbounded $A\subseteq \omega_1$, there is a
$\delta$ such that $A_\delta\subseteq A$.
\end{description}

\underline{Case 2.} $l=1$

For some stationary subset $S$ of $\omega_1\cap LIM$, there
is a sequence $\langle A_\delta:\,\delta\in S
\rangle$ such that
\begin{description}
\item{{\em (a)}} $A_\delta$ is an unbounded subset of $\delta$.
\item{{\em (b)}} $\otp{A_\delta}=\Upsilon$.
\item{{\em (c)}} For every unbounded $A\subseteq \omega_1$, there is a
$\delta$ such that $\card{A_\delta\setminus A}<\aleph_0$.
\end{description}

\underline{Case 3.} $l=2$

For some stationary $S\subseteq\omega_1\cap LIM$, there
is a sequence
\[
\langle \{A_\delta^n:\,n<\omega\}:\,\delta\in S
\rangle
\]
such that
\begin{description}
\item{{\em (a)}} Each $A^\delta_n$ is an unbounded subset of $\delta$.
\item{{\em (b)}} $\otp{A^\delta_n}=\Upsilon$.
\item{{\em (c)}} For every unbounded $A\subseteq \omega_1$, there is a
$\delta$ and an $n$ such that $A^\delta_n\subseteq A$.

\end{description}

\underline{Case 4.} $l=\bullet$.

For some stationary set $S\subseteq\omega_1\cap LIM$, there
is a sequence $\langle \{ A^\delta_m:\,m\le m^\ast(\delta)\}:\,\delta
\in S
\rangle$ such that
\begin{description}
\item{{\em (a)}} Each $A^\delta_m$ is an unbounded subset of $\delta$.
\item{{\em (b)}} $\otp{A^\delta_m}=\Upsilon$.
\item{{\em (c)}} For every unbounded $A\subseteq \omega_1$, there is a
$\delta$ and an $m\le m^\ast(\delta)$ such that $A^\delta_m\subseteq
A$.
\item{{\em (d)}} For all relevant $\delta$, we have
$m^\ast(\delta)<\omega$.
\end{description}

In the above, $LIM$ stands for the class of limit ordinals.

\end{Definition}

\begin{Remark}(1) One could, of course, consider the previous
definitions with $\omega_1$ replaced by some other uncountable
ordinal, in fact an uncountable regular cardinal. As our proofs only
deal with $\omega_1$, we only formulate our definitions in the
form given above.

Also, we could consider principles of the form
$\clubsuit^l_\Upsilon(T)$ in which
$T$ is a stationary subset of $\omega_1$ and parameter $\delta$ in the
above definitions is allowed to range only in $T$ (i.e. $S\cap T$).

(2) The definition that A. Ostaszewski [Ost] used for a
$\clubsuit$-sequence
$\langle A_\delta:\,\delta\in S\rangle$ requires that for each $A\in
[\omega_1]^{\aleph_1}$ there is a stationary set of $\delta$ such that
$A_\delta
\subseteq A$. It is well known that this is equivalent to our
definition
of $\clubsuit^0$. Hence 
$\clubsuit^0$ is the usual $\clubsuit$ principle of Ostaszewski, and we
shall
often
omit
the superscript $0$ when discussing this principle, and freely
use the equivalence between the definitions.

\end{Remark}

It is obvious that $\clubsuit^0_\Upsilon\implies\clubsuit^1_\Upsilon
\implies\clubsuit^2_\Upsilon$, and that
$\clubsuit^0_\Upsilon\implies\clubsuit^\bullet_\Upsilon\implies
\clubsuit^2_\Upsilon$. The result of the first sections \S2 and \S3
of the paper is that, except for the following
simple theorem, the above are the
only implications
that can be drawn. 

\begin{Theorem}\label{zfc}(1) Suppose that
$\Upsilon_1,\Upsilon_2<\omega_1$ are
limit ordinals and that $\clubsuit_{\Upsilon_1}$ and
$\clubsuit_{\Upsilon_2}$
both hold.

\underline{Then} $\clubsuit_{\Upsilon_1\cdot\Upsilon_2}$ holds.

(2)
$\clubsuit_{\Upsilon_1\cdot\Upsilon_2}\implies\clubsuit_{\Upsilon_1}$
for $\Upsilon_1$ limit $<\omega_1$ and $\Upsilon_2<\omega_1$. Similarly
for the other versions of $\clubsuit$ considered. 
\end{Theorem}

\begin{Proof}(1) Let $\langle A_\delta^l:\,\delta\in S_l\rangle$ for
$l=1,2$
exemplify
$\clubsuit_{\Upsilon_l}$. For $\delta\in \lim(S_1)\cap S_2$ 
we let
\[
B_\delta\deq\bigcup_{\alpha\in A^2_\delta}A^1_\alpha.
\]
Hence
$B_\delta$ is an unbounded subset of
$\delta$. 

Suppose that $A\in[\omega_1]^{\aleph_1}$.
For each $\alpha<\omega_1$, the set $A\setminus\alpha$
is an unbounded subset of $\omega_1$, hence contains stationarily
many $A^1_\delta$
as subsets.
So we can find an unbounded subset $T_1=T_1[A]$ of $S_1$ such that
\[
\alpha\in T_1\implies A^1_\alpha\subseteq
A\setminus\sup(T_1\cap\alpha).
\]
Now we can find a $\delta\in\lim(S_1)\cap S_2$ such that
$A^2_\delta\subseteq
T_1$.
Hence $B_\delta\subseteq A$ and $B_\delta$ is unbounded in $\delta$.
Moreover, $\otp{B_\delta}=\Upsilon_1\cdot\Upsilon_2$.

We have shown that $\langle B_\delta:\,\delta\in \lim(S_1)\cap
S_2\,\,\&
\,\,\otp{B_\delta}=\Upsilon_1\cdot\Upsilon_2\rangle$ witnesses
that $\clubsuit_{\Upsilon_1\cdot\Upsilon_2}$ holds (note that the
fact that the set of relevant $\delta$ is stationary follows
from the previous paragraph).

(2) Easy.
$\eop_{\ref{zfc}}$
\end{Proof}

The questions considered in the
paper are answered using the same basic technique,
with some changes in the definition of the particular forcing used. A detailed
explanation of the technique and the way it is used to prove that
$\club^1$ does not imply $\club^0$, is given in $\S\ref{vec}$.
The changes needed to obtain the other two theorems are presented
at the end of \S\ref{vec} and in
\S\ref{three}.

\section{Consistency of $\clubsuit^1$ and $\neg\clubsuit^0$}\label{vec}
\begin{Theorem}\label{1 not 0} $CON(\clubsuit^1+\neg\clubsuit)$.
\end{Theorem}

\begin{Proof} Throughout the proof, $\chi$ is a fixed large enough
regular
cardinal.

We start with a model $V$ of $ZFC$ such that
\[
V\models \diamondsuit(\omega_1) + 2^{\aleph_1}=\aleph_2,
\]
and use an iteration $\bar{Q}=\langle
P_\alpha,\name{Q}_\beta:\,\alpha\le\omega_2\,\,\&\,\,\beta<
\omega_2\rangle$.  The iteration is defined in the following
definition.

\begin{Definition}\label{order}
(1) By a candidate for a $\clubsuit$, we mean a sequence
of the form $\langle A_\delta:\,\delta<\omega_1\mbox{ limit }\rangle$,
such that
$A_\delta$ is an unbounded subset of $\delta$, with
$\otp{A_\delta}=\omega$.

(2) In $V$, we fix a continuously increasing sequence of countable
elementary
submodels of $(\HH(\chi),\in,
<^\ast_\chi)$, call it
$\bar{\bar{N}}=\langle N_i^0:\,i<\omega_1\rangle$,
such that $\HH(\aleph_1)\subseteq\bigcup_{i<\omega_1}N^0_i$ (this is
possible by
$CH$), and $\langle N^0_j:\,j\le i\rangle\in N^0_i$ for $i<\omega_1$.

(3) During the iteration, we do a bookkeeping which hands us
candidates for $\clubsuit$.

(4) Suppose that $\beta<\omega_2$, and let us define
$Q_\beta$, while
working in $V^{P_\beta}$.
\begin{enumerate}

\item Suppose that $CH$ holds in $V^{P_\beta}$ and the bookkeeping
gives us a
sequence $\bar{A}^\beta
=\langle A_\delta^\beta:\,\delta<\omega_1\mbox{  a limit ordinal
}\rangle$
which is a candidate for $\clubsuit$.
For some club $E_\beta$ of $\omega_1$
we choose a continuously increasing sequence $\bar{N}^\beta=\langle
N_i^\beta:\,i\in E_\beta\rangle$ of countable elementary
submodels of $(\HH(\chi),\in,<^\ast_\chi)$,
such that 
we have $\HH(\aleph_1)\subseteq\bigcup_{i\in E_\beta}N^\beta_i$,
and such that for every $i\in E_\beta$ we have $N^\beta_i\cap V=
N^0_i$, while $\langle N^\beta_j:\,j\le i\rangle
\in N^\beta_{\min(E_\beta\setminus (i+1))}$.
Furthermore, $\bar{A}^\beta\in N^\beta_{\min(E_\beta)}$.
{\em Then} $Q_\beta=Q_{\beta_{\bar{A},\bar{N}^\beta}}$
is defined by
\[
\begin{array}{llll}
Q_\beta\deq\{f: &(i)\, f \mbox{ is a partial function from
}\omega_1\mbox{ to
}\{0,1\}\\
&(ii)\,\otp{\Dom(f)}<\omega^\omega\\
&(iii)\,f\rest (N^\beta_i\cap\omega_1)\in
N^\beta_{\min(E_\beta\setminus (i+1))},
\mbox{ for }i\in E_\beta\\
&(iv)\,f^{-1}(\{1\})\cap A^\beta_\delta=\emptyset\implies
\card{\Dom(f)\cap A^\beta_\delta}<\aleph_0\\
&(v)\, f\in V\}
\end{array}
\]

\item If $\neg CH$, {\em then} $Q_\beta=\emptyset$. (Of course, 
our situation will be such that this case never occurs.) 

\end{enumerate}
In $Q_\alpha$, the order is given by
\[
f\le g\iff g\mbox{ extends }f\mbox{ as a function.}
\]

(5) For $\alpha\le\omega_2$, we define inductively 
\begin{eqnarray*}
P_\alpha\deq&\{p: &\,\Dom(p)\in[\alpha]^{\le\aleph_0}\,\,\&\,\,\left(
\forall\beta\in\Dom(p)\right)\\
&&(p(\beta)\mbox{ is a canonical hereditarily countable over Ord}\\
&&P_\beta\mbox{-name of
a member of }\name{Q}_\beta,\\
&&\mbox{ and }p\rest\beta\forces_{P_\beta}``p(\beta)\in\name{Q}_\beta")\}.
\end{eqnarray*}

The order in $P_\alpha$ is given by
\[
\begin{array}{llll}
p\le q\iff 
&(I)\,\,\Dom(p)\subseteq\Dom(q),\\
&(II)\,\,\mbox{For all }\beta\le\alpha,\mbox{ we have }
q\rest\beta\forces``p(\beta)\le q(\beta)",\\
&(III)\,\,
\{\gamma\in\Dom(p):\,p(\gamma)\neq q(\gamma)\}\mbox{ is finite }. 
\end{array}
\]
\end{Definition}

\begin{Definition} Suppose $\alpha\le\omega_2$, and $p\le q\in
P_\alpha$. {\em
Then}

(1) We say that $q$ {\em purely extends} $p$, if $q\rest\Dom(p)=p$. We
write
$p\le_{\pr}q$.

(2) We say that $q$ {\em apurely extends} $p$, if $\Dom(p)=\Dom(q)$. We
write
$p\le_{\apr}q$.

(3) The meaning of $p\ge_{\pr}q$ and  $p\ge_{\apr}q$ is defined in the
obvious
way.
\end{Definition}

\begin{Definition} Suppose that $\gamma<\omega_1$. A forcing notion $P$
is said
to be {\em purely} $\gamma$-{\em proper} if:

For every $p\in P$ and a continuously increasing sequence $\langle
N_i:\,i\le\gamma\rangle$ of countable elementary
submodels of $(\HH(\chi),\in,<^\ast_\chi)$ with $p,P\in
N_0$, $\langle N_j:\, j \le i\rangle\in N_{i+1}$, there is
a
$q\ge_{\pr} p$ which is $(N_i,P)$-generic for all $i\le\gamma$.
\end{Definition}

\begin{Fact}\label{trick} A $ccc$ forcing notion is purely
$\gamma$-proper
for every $\gamma<\omega_1$.
\end{Fact}
\begin{Proof of the Fact} This is because every condition in a $ccc$
forcing
is generic, see [Sh -f III, 2.6 and 2.9.]$\eop_{\ref{trick}}$
\end{Proof of the Fact}

{\em General facts about the iterations like the one we are using.}

\begin{Fact}\label{iteration} Iterations with the support we are using,
have the following general properties:

(1) $\alpha\le\beta\implies P_\alpha\subseteq P_\beta$ as ordered sets.

(2) $(\alpha\le\beta\,\,\&\,\,q\in P_\beta)\implies (q\rest \alpha\in
P_\alpha\,\,\&\,\,q\rest\alpha\le q).$

(3) $(\alpha\le\beta\,\,\&\,\,p\in P_\beta\,\,\&\,\,p\rest\alpha
\le q\in P_\alpha)\implies q\cup(p\rest [\alpha,\beta))\in P_\beta$ is
the least
upper bound of $p$ and $q$.

(4) If $\alpha<\beta$, {\em then} $P_\alpha\complB P_\beta$.
Hence, $G_{P_{\alpha+1}}/G_{P_\alpha}$ gives rise to
a directed subset of $Q_\alpha$ over $V[G_{P_\alpha}]$.

(5) If $\langle p_i:\,i<i^\ast<\omega_1\rangle$ is a
$\le_{\pr}$-increasing
sequence in $P_{\alpha^\ast}$ for some
$\alpha^\ast\le\omega_2$, {\em then}
$p\deq\bigcup_{i<i^\ast}p_i$ is a condition in $P_{\alpha^\ast}$
and for every $i<i^\ast$ we have $p_i\le_{\pr} p$.

(6) Pure properness is preserved by the iteration. 
Moreover, 
for any $\gamma<\omega_1$, pure $\gamma$-properness is preserved by the
iteration. 
\end{Fact}
\begin{Proof of the Fact} (1)-(5) Just checking.

(6) The statement follows from some more general facts proved in
[Sh -f, XIV]. A direct proof can be given along the
lines of the proof that countable support iterations preserve
properness,
[Sh -f, III 3.2]. 
$\eop_{\ref{iteration}}$
\end{Proof of the Fact}

{\em Back to our specific iteration.}

\begin{Claim}\label{ccc} Suppose $\alpha^\ast<\omega_2$. In
$V^{P_{\alpha^\ast}}$, the forcing $Q_{\alpha^\ast}$ has the $ccc$.
Moreover, it has the property of Knaster.
\end{Claim}

\begin{Proof of the Claim}
We fix such an $\alpha^\ast$ and work in $V^{P_{\alpha^\ast}}$.
We assume $CH$, as otherwise we have defined $Q_{\alpha^\ast}$ as an
empty set.

Hence 
sequences $
\bar{N}^{\alpha^\ast}\deq\langle N_i^{\alpha^\ast}:\,i\in
E_{\alpha^\ast}\rangle$
and
$\langle A^{\alpha^\ast}_\delta:\,\delta<\omega_1\mbox{ limit}
\rangle$ are given. Let
\[
E\deq\{\delta\in
E_{\alpha^\ast}:\,N_\delta^{\alpha_\ast}\cap\omega_1=\delta\},
\]
so $E$ is a club of $\omega_1$. Suppose that $q_\alpha\in
Q_{\alpha^\ast}$ for
$\alpha<\omega_1$ are given. Let
\[
A\deq\{\delta\in E:\,\mbox{for some }\alpha\in E
\setminus\delta
\mbox{ we have }\delta>\sup(\delta\cap\Dom(q_\alpha))\}.
\]
$A$ contains a final segment of $\acc(E)$, as otherwise
we can find an increasing sequence $\langle \delta_i:\,i<
\omega^\omega\rangle$ from $\acc(E)\setminus A$. Choose $\alpha
\ge\sup\{\delta_i:\,i<\omega^\omega\}$ with $\alpha\in E$.
Hence for all $i<\omega^\omega$ we have that $\delta_i=
\sup[\Dom(q_\alpha)\cap\delta_i]$, which is in contradiction with
$\otp{\Dom(q_\alpha)}<\omega^\omega$.

Let $C$ be a club such that $A\supseteq C$.
For $\delta\in C$, we fix an ordinal $\alpha_\delta$ witnessing that
$\delta\in
A$. So $\alpha_\delta\in E\setminus\delta$ and $\delta
>\sup(\delta\cap\Dom(q_{\alpha_\delta}))$. 

For $\delta\in C$, let $g(\delta)$ be defined as the minimal ordinal
$\in E$
such that $q_{\alpha_\delta}\in N^{\alpha_\ast}_{g(\delta)}$
(note that $g$ is well defined). Hence, the set of $\delta\in C$ which
are closed
under $g$, is a club of $\omega_1$. Call this club $C_1$.

Note that there is a stationary $S\subseteq C_1$ such that for some
$\xi^\ast$ we
have
\[
\delta\in S\implies\sup (\delta\cap\Dom(q_{\alpha_\delta}))=\xi^\ast.
\]
Now notice that for $\delta_1<\delta_2\in C_1$, we have
$\Dom(q_{\alpha_{\delta_1}})\subseteq
N_{\alpha_{\delta_2}}^{\alpha^\ast}\cap
\omega_1=\alpha_{\delta_2}.$ So, if $\delta_1<\delta_2\in S$,
we have
\[
\Dom(q_{\alpha_{\delta_1}})\cap \Dom(q_{\alpha_{\delta_2}})
\subseteq \alpha_{\delta_2}\cap  \Dom(q_{\alpha_{\delta_2}})
\subseteq
\xi^\ast.
\]
Now let $\delta^\ast\deq\min(S)$, so $\delta^\ast>\xi^\ast$. 
By {\em (iii)} in the definition of
$Q_{\alpha^\ast,\bar{N}^{\alpha^\ast}}$,
for every $\delta\in S$ we have
\[
q_{\alpha_\delta}\rest (\Dom(q_{\alpha_\delta})\cap\xi^\ast)
=\left(  q_{\alpha_\delta}\rest (\Dom(q_{\alpha_\delta})\cap
\delta^\ast)\right)\rest\xi^\ast\in N_{\min(E_{\alpha^\ast}
\setminus(\delta^\ast+1))}^{\alpha^\ast}.
\]
So, there are only countably many possibilities, hence we can find an
uncountable set of $\alpha_\delta$ such that $q_{\alpha_\delta}$ are
pairwise compatible.
$\eop_{\ref{ccc}}$
\end{Proof of the Claim}

\begin{Remark} $ccc$ orders like the one above were considered by
Abraham, Rubin
and Shelah in [ARSh 153].
\end{Remark}

\begin{Conclusion}\label{concll} For all $\alpha\le\omega_2$, the forcing
$P_\alpha$
is purely $\gamma$-proper for all $\gamma<\omega_1$.
\end{Conclusion}

[Why? By Fact \ref{trick}, Fact \ref{iteration}(6) and Claim
\ref{ccc}.]

\begin{Claim}\label{induction} The following hold for every
$\alpha^\ast\le\omega_2$:

(1) In $P_{\alpha^\ast}$, if $p\le r$, {\em then}
for some unique $q$ we have
\[
p\le_{\pr}q\le_{\apr}r\,\,\,\&\,\,\,(\alpha\in\Dom(q)\,\,\&\,\,q(
\alpha)\neq r(\alpha)\implies\alpha\in \Dom(p)).
\]

(2) The following is impossible in $P_{\alpha^\ast}$:

There is a sequence $\langle q_i:\,i<\omega_1\rangle$ which is
$\le_{\pr}$-increasing, but for which there is an antichain
$\langle r_i:\,i<\omega_1\rangle$ such that $q_i\le_{\apr}r_i
$.

(3)
If $p\in P_{\alpha^\ast}$ and $\name{\tau}$ is a $P_{\alpha^\ast}$-name
of an
ordinal, {\em then} there is $q\in P_{\alpha^\ast}$ with $p\le_{\pr}
q$, and a
countable antichain
$\name{I}\subseteq\{r:\,q\le_{\apr}r\}$ predense above $q$,
such that each $r\in \name{I}$ forces a value to $\name{\tau}$.

(4)  If $\alpha^\ast<\omega_2$, {\em then} $\forces_{P_{\alpha^\ast}}
``\card{\name{Q}_{\alpha^\ast}}=\aleph_1"$.

(5) If $\alpha^\ast<\omega_2$, {\em then}
$V^{P_{\alpha^\ast}}\satisfies CH$.

(6) $Q_{\alpha^\ast}$ is closed under finite unions of functions which
agree
on their common domain.

(7) $V^{P_{\alpha^\ast}}\satisfies 2^{\aleph_1}=\aleph_2$.

(8) $P_{\alpha^\ast}$ satisfies $\aleph_2$-cc.
\end{Claim}

\begin{Proof of the Claim} (1) Define $q$ by $q\deq p\cup (r\rest
(\Dom(r)\setminus \Dom(p))$. 

(2) We prove this by induction on $\alpha^\ast$. The case
$\alpha^\ast=0$ is vacuous, and 
if $\alpha^\ast$ is a successor ordinal, the statement easily follows
from the fact that each $Q_\alpha$ has the property of Knaster.

Suppose that $\alpha^\ast$ is a limit ordinal
and $\langle q_i:\,i<\omega_1\rangle$,
$\langle r_i:\,i<\omega_1\rangle$ exemplify a
contradiction to (2). For $i<\omega_1$ let
$w_i\deq\{\alpha\in \Dom(q_i):\,r_i(\alpha)\neq q_i(\alpha)\}$, hence
$w_i$ is a
finite set. Without loss of generality, we can assume that sets
$w_i\,(i<\omega_1)$ form a $\Delta$-system
with root $w^\ast$. Let $\beta^\ast\deq\Max(w^\ast)+1$, so
$\beta^\ast<\alpha^\ast$.

Now notice that 
\[
\alpha\in \Dom(r_i)\cap\Dom(r_j)\,
\,\,\&\,\,\neg(\forces_{P_{\alpha}}
``r_i(\alpha), r_j(\alpha)\mbox{ are compatible}")
\]
implies that $ \alpha\in w^\ast$,
for any $i,j<\omega_1$.
Hence, $\langle q_i\rest\beta^\ast:\,i<\omega_1\rangle$
and $\langle r_i\rest\beta^\ast:\,i<\omega_1\rangle$ exemplify that
(2)  fails at $\beta^\ast$, contradicting the induction hypothesis.

(3) We work in $V^{P_{\alpha^\ast}}$. Fix such $p$ and $\name{\tau}$.
Let $J$ be an antichain
predense above $p$, such that every $r\in J$ forces a value
to $\name{\tau}$.

We try to choose by induction on $i<\omega_1$ conditions $p_i,r_i$ such
that
\begin{itemize}
\item $p_0=p$,
\item $j<i\implies p_j\le_{\pr}p_i$,
\item $r_i\in J$,
\item $p_i\le_{\apr}r_i$,
\item $j<i\implies r_i\bot r_j$.
\end{itemize}

If we succeed,
(2) is violated, a contradiction.

So, we are stuck at some $i^\ast<\omega_1$. We can let $q\deq
p_{i^\ast}$ and
$I\deq\{r_i:\,i<i^\ast\}$.
 
(4) Obvious from the definition of $Q_{\alpha^\ast}$.

(5) Can be proved by induction on $\alpha^\ast$, using
(3) and (4).

(6) Just check.

(7) Follows from the definition of $P_{\alpha^\ast}$, part (3)
of this claim, and
the fact that $V\satisfies 2^{\aleph_1}=\aleph_2$.

(8) By \ref{order}(5) and part (4) of this claim (see [Sh -f], III 4.1
for the analogue in the case of countable support iterations).
$\eop_{\ref{induction}}$
\end{Proof of the Claim}

\begin{Claim}\label{book} It is possible to arrange the bookkeeping, so
that
$\forces_{P_{\omega_2}}\neg\clubsuit.$
\end{Claim}

\begin{Proof of the Claim} As usual, using Claim \ref{induction}(7),
it suffices to prove that for
every $\alpha^\ast<\omega_2$, in $V^{P_{\alpha^\ast}}$
we have
\[
\forces_{Q_{\alpha^\ast}}``\langle
A^{\alpha^\ast}_\delta:\,\delta<\omega_1\rangle
\mbox{ is not a }\clubsuit\mbox{-sequence}."
\]
Let $G$ be $Q_{\alpha^\ast}$-generic
over $V^{P_{\alpha^\ast}}$, and let $F\deq\bigcup G$. Let
$A\deq F^{-1}(\{0\})$. Suppose that $A\supseteq A^{\alpha^\ast}_\delta$
for some
$\delta$. Then for every $f\in G$
we have $f^{-1}(\{1\})\cap A^{\alpha^\ast}_\delta=\emptyset$,
so $\card{\Dom(f)\cap A^{\alpha^\ast}_\delta}<\aleph_0$.

However, the following is true:

\begin{Subclaim}\label{density} The set 
\[
\II\deq\{f\in Q_{\alpha^\ast}:\,\card{\Dom(f)\cap
A^{\alpha^\ast}_\delta}=\aleph_0
\mbox{ or }f^{-1}(\{1\})\cap A^{\alpha^\ast}_\delta\neq\emptyset\}
\]
is dense in $Q_{\alpha^\ast}$.
\end{Subclaim}
\begin{Proof of the Subclaim} Given $f\in Q_{\alpha^\ast}$.
If
$\Dom(f)\cap A^{\alpha^\ast}_\delta$ is infinite, then $f\in \II$.
Otherwise,
let $\beta\deq\min(A^{\alpha^\ast}_\delta)\setminus\Dom(f)$. Let $g\deq
f\bigcup\{(\beta,1)\}$, hence $g\ge f$ and $g\in \II$.
$\eop_{\ref{density}}$
\end{Proof of the Subclaim}  
We obtain a contradiction, hence $A$ is not a superset of
$A^{\alpha^\ast}_\delta$.
$\eop_{\ref{book}}$
\end{Proof of the Claim}

\begin{Definition}\label{bad}
 Suppose that
\begin{description}
\item{(a)} $\gamma<\omega_1$,
\item{(b)} $\bar{N}=\langle N_i:\,i\le\gamma\rangle$ is a continuous
increasing
sequence of countable elementary submodels of $\langle
\HH(\chi),\in,<^\ast_\chi\rangle$,
\item{(c)} $\name{\tau},\bar{Q}\in N_0$ and $p\in P_{\omega_2}\cap
N_0$, 
\item{(d)} $p\forces``\name{\tau}\in [\omega_1]^{\aleph_1}"$ and
\item{(e)} $\bar{N}\rest(i+1)\in N_{i+1}$ for $i<\gamma$.
\end{description}

We say that $\varepsilon\le\gamma$ is {\em bad}
for $(\bar{N},\name{\tau},p,\bar{Q})$ if $\varepsilon$ is a limit ordinal,
and there
are no $r_n$, $\beta_n\in N_\varepsilon\,(n<\omega)$ such that 
\begin{description}
\item{(1)} $ r_n\forces_{P_{\omega_2}}``\beta_n\in \name{\tau}$",
\item{(2)} $\bigcup_{n\in \omega}\beta_n=N_\varepsilon\cap\omega_1,$
\item{(3)} $ r_n\ge p \mbox{ for all } n,$
\item{(4)} $\beta_n \mbox{ increase with }n,$ 
\item{(5)} $\mbox{ for some }n_0\in\omega$
the set $\{r_n:\,n\ge n_0\}$
has an upper bound in $P_{\omega_2}$
\item{(6)} ${\bar{r}}_{\bar{N}\rest\varepsilon,p,\name{\tau}}
\deq\langle r_n:\,n<\omega\rangle$ and 
${\bar{\beta}}_{\bar{N}
\rest\varepsilon,p,\name{\tau}}\deq\langle \beta_n:\,n<\omega\rangle$
 are definable in $(\HH(\chi)^V,\in,<^\ast_\chi)$ from the isomorphism
type of
$(\langle N_\xi:\,\xi\le\varepsilon\rangle,
 p,\name{\tau},\bar{Q})$ (we shall sometimes
abbreviate this by saying that these objects are defined in a
{\em canonical} way).
\end{description}
\end{Definition}

\begin{Main Claim}\label{2.5} Suppose that $\bar{N},\gamma,p$ and
$\name{\tau}$
are as in Definition \ref{bad}.
\underline{Then} the set
\[
B\deq\{\varepsilon\le\gamma:\,\varepsilon\mbox{ bad for
}(\bar{N},\name{\tau},p,\bar{Q})\}
\]
has order type $<\omega^\omega$.
\end{Main Claim}

\begin{Proof of the Main Claim} We start by 

\begin{Subclaim}\label{p} Let $\bar{N},\gamma, p$ and $\name{\tau}$ be
as in the
hypothesis of Claim \ref{2.5}. \underline{Then}, 
we can choose canonically a sequence $\bar{p}=\langle
p_j:\,j<\omega\gamma\rangle$ such that
\begin{enumerate}
\item $\bar{p}$ is $\le_{\rm pr}$-increasing,
\item $p_0=p$,
\item\label{subset} For $i<\gamma$ and $n<\omega$,
we have that $p_{\omega i+n}\in
N_{i+1}$.
\item\label{glavni}  For each $i<\gamma$,
for every formula $\psi(x,y)$ with parameters
in $N_{i}$, there
are infinitely many $n$  such that one of the following occurs: 

\begin{description}
\item{$(\alpha)$} For no $p'\ge p_{\omega i +n}$ do we have that for
some $y$,
the formula $\psi(p',y)$ holds. 

\item{$(\beta)$} For the $<^\ast_\chi$-first $r\ge p_{\omega i+n}$ such
that
$\psi(r,y)$ holds
for some $y$, we have $r\ge_{\apr} p_{\omega i+n+1}$.

\end{description}
\item\label{limits} For $j<\omega\gamma$ a limit ordinal, we have
$p_j=\cup_{i<j}
p_i$.
\end{enumerate}
\end{Subclaim}

\begin{Proof of the Subclaim} We prove this by induction on $\gamma$, 
for all
$\bar{N}$ and $p$.

If $\underline{\gamma=0}$, there is nothing to prove.

If $\underline{\gamma<\omega_1}$ is a limit ordinal, we fix an
increasing 
sequence $\langle\gamma_k:\,k<\omega\rangle$ which is cofinal
in $\gamma$, such that $\gamma_0=0$
(we are taking the $<_\chi^\ast$-first sequence like that).
By induction on $k$ we define $\langle
p_j:\,\omega\gamma_k< j\le\omega{\gamma_{k+1}}\rangle$. We let $p_0\deq
p$.
At the stage $k$ of the induction we 
use the induction hypothesis
with $p_{\omega\gamma_k}, \langle N_j:\,\omega\gamma_k<j\le
\omega\gamma_{k+1}\rangle$ here standing for $p,\bar{N}$ there,
obtaining $\langle p_j:\,\omega\gamma_k<j\le\omega\gamma_{k+1}
\rangle$, noticing that $p_{\omega\gamma_k}\in N_{\omega\gamma_{k+1}}$. 
We define $p_{\omega\gamma_{k+1}}\deq\bigcup_{j<\omega\gamma_k}p_j$.
We thus obtain
$\langle p_j:\,\omega\gamma_k<j\le\omega\gamma_{k+1}\rangle$ in $V$.
As the parameters used
are in $N_{\omega\gamma_k +1}$, by the fact that our choice is
canonical, we
have that $\langle
p_j:\,\omega\gamma_k<j\le\omega\gamma_{k+1}\rangle\in
N_{\omega\gamma_{k+1}+1}$.

Suppose that $\underline{\gamma=\gamma'+1}$. 
By the induction hypothesis, we can find
a sequence $\langle p_j:\,j<\omega\gamma'\rangle$
satisfying the subclaim for  $p$ and $\bar{N}\rest\gamma'$.
As $\bar{N}\rest\gamma\in N_\gamma$, again we have that
the sequence $\langle p_j:\,j<\omega\gamma'\rangle$ is in $N_\gamma$.
Let
$p_{\omega\gamma'}\deq \cup_{j<\omega\gamma'} p_j$.

We list as $\langle \psi^{\gamma}_n=\psi_n:\,n<\omega\rangle$ all
formulas
$\psi(x,y)$ with
parameters in $N_{\gamma'}$,
so that each formula appears infinitely often,
picking the $<_\chi^\ast$-first such enumeration. By induction on
$n<\omega$, we
choose $p_{\omega\gamma'+n}$. We have already chosen
$p_{\omega\gamma'}$.

At the stage $\underline{n+1}$ of the induction, we consider $\psi_n$.
If $(\alpha)$ holds, we just let
$p_{\omega\gamma'+n+1}\deq p_{\omega\gamma'+n}$.
Otherwise, there is a condition $r\ge p_{\omega\gamma'+n}$
such that  $\psi_n(r,y)$ for some $y$. By elementarity,
the $<^\ast_\chi$-first such $r$ is in $N_{\gamma'+1}$.
By Claim \ref{induction}(1),  there is a unique $q$ such that
$r\ge_{\apr}q\ge_{\pr}p_{\omega\gamma'+n}$ and
$\alpha\in\Dom(q)\,\,\&\,\,r(\alpha)\neq q(\alpha)\implies
\alpha\in\Dom(p)$. Hence, $q\in
N_{\gamma'+1}$ and we set 
$p_{\omega\gamma+n+1}\deq q$. $\eop_{\ref{p}}$
\end{Proof of the Subclaim}

We now choose $\bar{p}$ as in the Subclaim, using our fixed
$\gamma,\bar{N},
\name{\tau}$ and $p$. 

\begin{Note}\label{after} For every limit $\varepsilon<\gamma$ we have
that
$\Dom(p_{\omega\varepsilon})=N_{\varepsilon}\cap\omega_2$.
\end{Note}

[Why? Let $i<\omega\varepsilon$ be given, and let $\alpha\in
N_{i}\cap\omega_2$.
Consider the formula
$\psi(x,y)$ which says that $x=y\in P_{\omega_2}$ and
$\alpha\in \Dom(x)$. This is a formula with parameters in $N_{i}$.
Option $(\alpha)$ from item
\ref{glavni}. of Subclaim \ref{p} does not occur, so there is $m$ and
$r\ge_{\apr}p_{\omega i +m}$ such that $\psi(r,y)$ holds for some $y$.
Hence
$\alpha\in\Dom(r)=\Dom(p_{\omega i +m})\subseteq
\Dom(p_{\omega (i+1)})$. So $N_i\cap\omega_2\subseteq
\Dom(p_{\omega(i+1)})$,
and hence $N_\varepsilon\cap\omega_2\subseteq\Dom(p_{\omega\varepsilon})$.

On the other hand,
if $\alpha\in \Dom(p_{\omega\varepsilon})$, there is $i<\varepsilon$ such
that
$\alpha\in\Dom(p_{\omega i})\subseteq N_{i+1}\subseteq N_{\varepsilon}$.]

\begin{Observation}\label{inV} Suppose $\alpha\le\omega_2$,
while $q\in P_\alpha$ and $w\in
[\Dom(q)]^{<\aleph_0}$.
\underline{Then} there is $q^+\ge q$ in $P_\alpha$ such that
\begin{description}
\item{$(\ast)^\alpha$} If $i\in w\cup\{j\in \Dom(q)
:\,q(j)\neq q^+(j)\}$, then $q^+(i)\in V$ (an object),
and not just $q^+\rest i\forces ``q^+(i)\in V"$ (not just a name).
\end{description}
\end{Observation}

[Why? By induction on $\alpha$. The induction is
trivial for $\alpha=0$, and in the case of $\alpha$ a limit ordinal it follows
from the finiteness of $w$. Suppose that $\alpha=\beta+1$. We have
$q\rest\beta\forces ``q(\beta)\in V"$, so we can find $r\in P_\beta$ such that
$r\ge q\rest\beta$, and $A$ such that $r\forces``q(\beta)=A"$. Now apply
$(\ast)^\beta$ with $r$ in place of $q$ and $(w\cap\beta)\cup\{j:\,
r(j)\neq q(j)\}$ to obtain $q^+_\beta$. Let $q^+\deq q^+_\beta
\frown\{\langle \beta,A\rangle\}$.]

\medskip

{\em Continuation of the proof of \ref{2.5}.}
 
\medskip

Since $\bar{p}$ is $\le_{\pr}$-increasing, the limit of $\bar{p}$ is a
condition,
say $p_\ast$. 
Now let $q^\ast\ge p_\ast$ be the $<^\ast_\chi$-first such that
$q^\ast\forces
``\beta\in\name{\tau}"$ for some $\beta>N_\gamma\cap\omega_1$, and with the
property
\[
[\alpha\in \Dom(p_\ast)\,\,\&\,\,p_\ast(\alpha)\neq q^\ast(\alpha)]\implies
q^\ast(\alpha)\mbox{ an object},
\]
which exists by Observation \ref{inV}.
Let
$w^\ast \deq\{\alpha\in \Dom(p_\ast):\,p_\ast(\alpha)\neq
q^\ast(\alpha)\}$.


\medskip

We now define
\[
b\deq\{\varepsilon\le\gamma:\,\left(\bigcup_{\alpha\in
w^\ast}\Dom(q^{\ast}(\alpha))\cap (N_\varepsilon\cap\omega_1)\right)\mbox{
is
unbounded in }N_\varepsilon
\cap\omega_1\}.
\]

\begin{Note} ${\rm otp}(b)<\omega^\omega$.
\end{Note}

[Why? Suppose that $\varepsilon_j$ for $j<\omega^\omega$ are elements of
$b$,
increasing with $j$. Now, for every $j<\omega^\omega$ we know that
$N_{\varepsilon_j}\cap\omega_1$ is bounded in
$N_{\varepsilon_{j+1}}\cap\omega_1$, but
$(\bigcup_{\alpha\in w^\ast}\Dom(q^{\ast}(\alpha))\cap
(N_{\varepsilon_{j+1}}\cap\omega_1)\mbox{ is unbounded in
}N_{\varepsilon_{j+1}}
\cap\omega_1$. Hence $\bigcup_{\alpha\in
w^\ast}\Dom(q^{\ast}(\alpha))\cap 
[N_{\varepsilon_j}\cap\omega_1,N_{\varepsilon_{j+1}}\cap\omega_1)
\neq\emptyset$. However,
by the definition of the forcing, $\otp{\bigcup_{\alpha\in
w^\ast}\Dom(q^{\ast}(\alpha))}
<\omega^\omega$, a contradiction.]

\medskip

{\em Continuation of the proof of \ref{2.5}.}

\medskip

Our aim is to show that  $B\subseteq b$
($B$ was defined in the statement of the Main Claim). So, let
$\varepsilon^\ast\in(\gamma+1)\setminus b$ be a limit ordinal.
We show that $\varepsilon^\ast\notin B$. We have to
define $\bar{r}\deq\bar{r}_{\bar{N}\rest{\varepsilon^\ast},p,\name{\tau}}$
and
$\bar{\beta}\deq\bar{\beta}_{\bar{N}\rest{\varepsilon^\ast},p,\name{\tau}}$
so to
satisfy (1)--(5) from the definition of $B$, and to do so in a
canonical way,
to be able to prove Subclaim \ref{X} below, hence showing
that (6) from Definition \ref{bad} holds.

Let
$\xi\deq[\sup\left(\bigcup_{\alpha\in
w^\ast}\Dom(q^{\ast}(\alpha)\right)
\cap N_{\varepsilon^\ast}\cap\omega_1]+1$, so $\xi<N_{\varepsilon^\ast}\cap
\omega_1$. We
enumerate $N_{\varepsilon^\ast}\cap w^\ast$ as
$\{\alpha_0,\ldots,\alpha_{n^{\ast}-1}\}$. 
By Note \ref{after}, we can fix $j^\ast<\varepsilon^\ast$
such that $\{\alpha_0,\ldots,\alpha_{n^\ast-1}\}$
$\subseteq\Dom(p_{\omega j^\ast})
$. Let $j^\ast$ be the first such.
Also let $\delta\deq N_{\varepsilon^\ast}\cap\omega_1$.
Now we observe that for all $l<n^\ast$, we have $q^\ast(\alpha_l)\rest\xi\in 
N_{\varepsilon^\ast}$.

[Why? Clearly, there is $\varepsilon'<\varepsilon^\ast$ such that
$\{\alpha_0,\ldots,\alpha_{n^\ast-1},\xi\}\subseteq N_{\varepsilon'}$.
With $\bar{\bar{N}}$ defined in Definition \ref{order}(2), we have that
$\bar{\bar{N}}\in N_0$. Also, we have that 
\[
\emptyset
\forces_{\alpha_{n^\ast-1}}``\name{E}\deq\bigcap_{l<n^\ast}\name{E}_{\alpha_l}
\mbox{ is a club of }\omega_1",
\]
(cf. Definition \ref{order}(4)1). Hence, by properness and the choice
of $\bar{N}$, we have that for every $\varepsilon\in [\varepsilon',
\gamma]$, we have that
\[
\emptyset
\forces_{\alpha_{n^\ast-1}}``N_\varepsilon\cap \omega_1\in \name{E}".
\]
Let $i\deq N_{\varepsilon'}\cap\omega_1$, hence $N^0_i\in N_{\varepsilon'+1}$.
In particular, we have $\emptyset
\forces_{\alpha_{n^\ast-1}}``i\in \name{E}"$ and $N^0_i\cap\omega_1
<N_{\varepsilon^\ast}\cap\omega_1$. So for all $l<n^\ast$ we have
$q^\ast(\alpha_l)\rest\xi=q^\ast(\alpha_l)\rest(N^0_i\cap\omega_1)$, but
\[
\emptyset
\forces_{\alpha_l}``N^0_i\cap\omega_1=\name{N}^{\alpha_l}_i\cap\omega_1",
\]
hence by Definition \ref{order}(4)1.$(iii)$, we have
$q^\ast\rest\alpha_l\forces``q^\ast(\alpha_l)\rest\xi\in
N^0_{\min(\name{E}_{\alpha_l}\setminus(i+1))}$. But
$\emptyset_{\alpha_l}\forces``\min(\name{E}_{\alpha_l}\setminus(i+1))\in
N_{\varepsilon'+1}[\name{G}]"$, hence
$q\rest\alpha_l\forces``q(\alpha_l)\rest\xi\in N_{\varepsilon'+1}[\name{G}]"$.
By properness and the fact that $q^\ast(\alpha_l)\in V$, we have
$q^\ast(\alpha_l)\rest\xi\in N_{\varepsilon'+1}$.]

Let us pick the $<^\ast_\chi$-first increasing sequence
$\langle\varepsilon_n:\,n<\omega\rangle$ such that
$\varepsilon^\ast=\bigcup_{n<\omega}\varepsilon_n$, while
$\omega
j^\ast+1<\varepsilon_0$ and $\xi\in N_{\varepsilon_0}$, in addition to
$(\forall l<n^\ast)[q^\ast(\alpha_l)\rest\xi \in N_{\varepsilon_0}]$. 
\smallskip

\underline{Defining $r_n$ and $\beta_n$}.
We do this by induction on $n$.
If $\underline{n=0}$, we set $r_0\deq p_{\omega\varepsilon_0}$, and also
let
$m_0=0$, $\xi_0=\xi$.

At stage $\underline{n+1}$, we assume that at the stage $n$
we have chosen $r_n\in N_{\varepsilon_n+1}\cap
P_{\omega_2}$ and $m_n<\omega$ so that $r_{n}\ge_{\apr}
p_{\omega\varepsilon_n
+m_n}$.
We also have chosen $\xi_n,\beta_n\in N_{\varepsilon_n+1}$.

We define a formula $\varphi_n(x,y)$ which says

\begin{enumerate}
\item\label{prvi} $x\in P_{\omega_2}$ and $y$
is an ordinal $>\Max\{\beta_n,N_{\varepsilon_n}\cap\omega_1\}$.
\item\label{drugi} $x\forces``y'\in\name{\tau}"$ for some $y'>y$.
\item\label{peti} If $l<n^\ast$, then $x(\alpha_l)$ is an object, not
a name, and $x(\alpha_l)\rest\xi=q^\ast(\alpha_l)\rest\xi$.
\item\label{treci} For $l<n^\ast$, we have
$x(\alpha_l)\rest\xi\in N_{\varepsilon_0}$ and
$\Dom(x(\alpha_l))\setminus\xi\subseteq\omega_1\setminus\xi_n$.
\item\label{cetvrti}
For all $\alpha$ we have
\begin{eqnarray*}
\alpha\in\Dom(x)\cap\Dom(p_{\omega\varepsilon_n + m_n})
&\&\,\,x(\alpha)\neq  p_{\omega\varepsilon_n + m_n}(\alpha)\\
&\implies\alpha\in\{\alpha_0,\ldots\alpha_{n^\ast-1}\}.
\end{eqnarray*}

\end{enumerate}

Hence, $\varphi_n$ is a formula with parameters in
$N_{\varepsilon_n+1}\subseteq
N_{\varepsilon_{n+1}}$. Also,
we have that
$\varphi_n(q^\ast,\delta)$ holds. 

By the choice of $\bar{p}$, there is $m_{n+1}>m_n$
(we pick the first one) such that for the $<^\ast_\chi$-first $r
\ge p_{\omega(\varepsilon_{n+1})+ m_{n+1}-1}$ for which there is $y$ for
which
$\varphi_n(r,y)$ holds, we have
$r\ge_{\apr} p_{\omega(\varepsilon_{n+1})+ m_{n+1}}$. We let
\[
r_{n+1}
\deq r\cup
(p_{\omega\varepsilon_{n+1}+m_{n+1}}
\rest \Dom(p_{\omega\varepsilon_{n+1}+m_{n+1}})
\setminus \Dom (r)).
\]
Note that $r_{n+1}\in N_{\varepsilon_{n+1}+1}$ and that
$\varphi_n(r_{n+1},y)$ must hold for some $y$.
The $<^\ast_\chi$-first such $y$
is an element of $N_{\varepsilon_{n+1}+1}$, and we choose it to be
$\beta_{n+1}$. 

Finally, we define $\xi_{n+1}\deq\min\left( N_{\varepsilon_{n+1}}\setminus
\sup\{\bigcup_{l<n^\ast}\Dom(r_{n+1}(\alpha_l))\setminus\xi\}\right)$.

\medskip

At the end, we obtain (canonically chosen) sequences $\langle
r_n:\,n<\omega\rangle$, $\langle \beta_n:\,n<\omega\rangle$,
$\langle \xi_n:\,n<\omega\rangle$ and $\langle m_n:\,n<\omega\rangle$
such that

\begin{enumerate}
\item\label{uno} $r_n\ge_{\apr} p_{\omega\varepsilon_n+m_n}$.
\item\label{due} $\xi_0=\xi$ and $\xi_n$ are strictly increasing with
$n$.
\item\label{tre} For all $l<n^\ast$, we have
$\Dom(r_n(\alpha_l))\setminus\xi
\subseteq (\xi_n,\xi_{n+1})$ and $r_n(\alpha_l)$ is an object.
\item\label{cete} $r_n\forces_{P_{\omega_2}}``\beta_n\in\name{\tau}"$.
\item\label{cinkve} $\beta_{n+1}>\beta_n$.
\item\label{se} $\bigcup_{n<\omega}\beta_n=
N_{\varepsilon^\ast}\cap\omega_1$.
\item\label{sete} $r_n\in N_{\varepsilon^\ast}$.
\item\label{oce} For $l<n^\ast$, we have $r_n(\alpha_l)\rest\xi=
r_1(\alpha_l)\rest\xi$.
\item\label{nonne} $\alpha\in\{\beta\in \Dom(r_n):\,r_n(\beta)\neq
p_{\omega\varepsilon_n+m_n}(\beta)\}\implies \alpha\in
\{\alpha_0,\ldots\alpha_{n^\ast}-1\}$.
\end{enumerate}

[Why?
By item \ref{cetvrti}. in the
definition of $\varphi_n$.]

We will use $r_n,\beta_n\,(n<\omega)$ to witness that
$\varepsilon^\ast\notin B$. It is true that $r_n\ge p$
and $\beta_n$  increase with $n$, and their limit is
$N_{\varepsilon^\ast}\cap
\omega_1$. 
We need to show that for some $n_0$, the
sequence $r_n\,(n\ge n_0)$ has an upper bound in $P_{\omega_2}$.
The natural choice to use would be $\bigcup_{n<\omega}r_n$, but
this is not necessarily a condition! 

[Why? By item \ref{nonne}. above, all $r_n$ for $n>0$
agree on $\alpha$ such that
$\alpha\notin\{\alpha_0,\ldots,\alpha_{n^{\ast}-1}\}$. 
By items \ref{due}, \ref{tre}. and \ref{oce} above, we even know that
for every $l<n^{\ast}$, the union $\bigcup_{n<\omega}r_n(\alpha_l)$ is
a
function.
If $\delta'<N_{\varepsilon^\ast}\cap\omega_1$,  then
for all $l<n^\ast$ we have
$\bigcup_{n<\omega}r_n(\alpha_l)\rest\delta'
=\bigcup_{n<n'}r_n(\alpha_l)\rest\delta'$ for some $n'<\omega$,
so this is a condition in $Q_{\alpha_l}$ (by Claim \ref{induction}
(6)). If
$\delta'>N_{\varepsilon^\ast\cap\omega_1}$, then
$\bigcup_{n<\omega}r_n(\alpha_l)\cap\delta'$ is finite. However,
it is possible that for some $\alpha_l$ 
it is forced that the intersection of the set
$\bigcup_{n\in\omega}\Dom(r_n(\alpha_l))$
with $\name{A}^{\alpha_l}_{N_{\varepsilon^\ast\cap\omega_1}}$ is infinite,
so $\bigcup_{n<\omega}r_n(\alpha_l)$ might fail to be a condition in
$\name{Q}_{\alpha_l}$.]

(We remark that
it is because of this point that we are getting $\clubsuit^1$ and not
$\clubsuit$
in $V^P$.)

Now, we define conditions $q_l^\ast$ for $l\le n^{\ast}$ as follows.
First set $\alpha_{n^{\ast}}\deq \omega_2$. By
induction on
$l\le n^{\ast}$
we choose $q^\ast_l\in P_{\alpha_l}$, so that
\begin{description}
\item{$(a)$} $q_l^\ast\le q^\ast_{l+1}$,
\item{$(b)$} $q_l^\ast\rest\alpha_l$
is above $r_n\rest \alpha_l$ for all $n$ large
enough.
\end{description}
This clearly suffices, as $q_{n^{\ast}}^\ast\cup q^\ast\rest
(\Dom(q^\ast)\setminus \Dom(q_{n^{\ast}}^\ast))$ is a condition
in $P_{\omega_2}$ which is above all but finitely many $r_n$.

{\em The choice of $q_l^{\ast}$.} Let $q^\ast_0\deq q^{\ast}\rest
\alpha_0=
p_\ast\rest \alpha_0$.
Given $q_l^\ast\in P_{\alpha_l}$, with $l<n^{\ast}$. We can find
$q_l^{\ast\ast}\ge q_l^\ast$ in $P_{\alpha_l}$,
such that $q_l^{\ast\ast}\forces``\min
(\name{A}^{\alpha_l}_{N_{\varepsilon^\ast\cap\omega_1}}\setminus
\xi_0)=\zeta_l$"
for some ordinal $\zeta_l$.
By item \ref{tre}. above, the ordinal $\zeta_l$ belongs to
$\Dom(r_n(\alpha_l))$
for at most one $n$. Let $n_l$ be greater than this $n$.
Hence there is a condition $q_l^+$ in $P_{\alpha_l+1}$ such that
$q_l^+(\alpha_l)$ is an object and
\[
q_l^+\rest
\alpha_l=q_l^{\ast\ast}\,\,\,\&\,\,\,q_l^+(\alpha_l)\ge\bigcup_{n\ge
n_l}r_n(\alpha_l)\,\,\,\&\,\,\,q_l^+(\alpha_l)(\zeta_l)=1.
\]
Now let $q^\ast_{l+1}\deq q_l^+\cup\bigcup_{n\ge
n_l}r_n\rest[\alpha_l+1,
\alpha_{l+1})$. 
Note that $q_{l+1}^\ast(\alpha)$ is forced to be
a function, for any $\alpha\in
\Dom(q_l)$, as
all $r_n$ agree
on $[\alpha_l+1,\alpha_{l+1})$. Also, $q^\ast_{l+1}(\alpha)$
is forced to be in $V$.

Now, the sequence $\langle q^\ast_{\alpha_l}:\,l\le n^{\ast}\rangle$ is
as
required.

To finish the proof of the Main Claim, we need to observe

\begin{Subclaim}\label{X} Suppose that $\bar{N}$ and $\bar{M}$
are two equally long countable continuously increasing sequences of
countable elementary submodels of $\langle
\HH(\chi),\in,<^\ast_\chi,p,\name{\tau},\bar{Q}\rangle$
with 
$\bar{Q}^N=\bar{Q}^M=\bar{Q}$, and $F=\langle f_i:\,i<{\rm lg}(\bar{N})\rangle$
is an increasing sequence of isomorphisms $f_i:\,N_i\into M_i$.

\underline{Then}, if $\bar{\beta}_{\bar{N},p,\name{\tau}}$ and
$\bar{r}_{
\bar{N},p,\name{\tau}}$ are defined, so are 
$\bar{\beta}_{\bar{M},F(p),F(\name{\tau})}$ and
$\bar{r}_{\bar{M},F(p),F(\name{\tau})}$. 
Moreover,
$\bar{\beta}_{\bar{M},F(p),F(\name{\tau})}=\bar{\beta}_{\bar{N},p,\name{\tau}}$
and $\bar{r}_{\bar{M},F(p),F(\name{\tau})}=
F(\bar{r}_{\bar{N},p,\name{\tau}})$.
\end{Subclaim}
\begin{Proof of the Subclaim} Check, looking at the way
$\bar{\beta},\bar{r}$ were defined.$\eop_{\ref{X}}$
\end{Proof of the Subclaim}
 $\eop_{\ref{2.5}}$
\end{Proof of the Main Claim}

To finish the proof of the Theorem, we prove

\begin{Claim}\label{conc} $\forces_{P_{\omega_2}}\clubsuit^1$.
\end{Claim}

\begin{Proof of the Claim}
We use the following equivalent reformulation of
$\diamondsuit$ in $V$: 

There is a sequence
\[
\left\langle \bar{N}^\delta=\langle N^\delta_i:\,i<\delta\rangle:
\,\delta<\omega_1\right\rangle,
\]
such that
\begin{enumerate}
\item Each $\bar{N}^\delta=\langle N^\delta_i:\,i<\delta\rangle$
is a continuously increasing sequence of countable elementary
submodels of $\langle
\HH(\chi),\in,<^\ast_\chi,p,\name{\tau},\bar{Q}\rangle$,
with $N_i^\delta\cap\omega_1<\delta$ and $\bar{N}^\delta\rest(i+1)
\in N^\delta_{i+1}$ for $i<\delta$.
Here, $p,\bar{Q}$ and $\name{\tau}$ are constant symbols.
In addition, $\bar{Q}^{N_0^\delta}=\bar{Q}$.

\item For every continuously increasing sequence
$\bar{N}=\langle N_i:\,i<\omega_1\rangle$
of countable elementary submodels
of $\langle 
\HH(\chi),\in,<^\ast_\chi,p,\name{\tau},\bar{Q}\rangle$
such that $\bar{Q}^{N_0}=\bar{Q}$,
there is a stationary set of $\delta$ such that
for all $i<\delta$ 
the isomorphism type of $N_i$ and $N_i^\delta$ is the same,
as is witnessed by some sequence of isomorphisms
$\langle f^\delta_i:\,i<\delta\rangle$ which is increasing with $i$.
\end{enumerate}

For each limit ordinal $\delta$, let
$N^\delta\deq\bigcup_{i<\delta}N_i^\delta$.
We define $A_\delta$:

If $\bar{\beta}_{\bar{N}^\delta,p^{N_\delta},\name{\tau}^{N_\delta}}$
is well
defined, then we let $A_\delta\deq
{\rm{Rang}}(\bar{\beta}_{\bar{N}^\delta,p^{N_\delta},\name{\tau}^{N_\delta}})$.
Otherwise, we let
$A_\delta$ be the range of any cofinal $\omega$-sequence in $\delta$.
Note that
in any case $A_\delta$ is an unbounded subset of $\delta$
of order type $\omega$.

We claim that $\langle A_\delta:\,\delta<\omega_1\rangle$
exemplifies that $V^P\models\clubsuit^1(\omega_1)$.We have to check
that for every unbounded subset $A$ of $\omega_1$ in
$V^{P_{\omega_2}}$, there is
a $\delta<\omega_1$ with
$\card{A_\delta\setminus A}<\aleph_0$.

Suppose this is not true. So, there are $p^\ast,\name{\tau}^\ast$
exemplifying
this, that is
\[
p^\ast\forces``\name{\tau}^\ast\in [\omega_1]^{\aleph_1}\mbox{ and for
all
}\delta\mbox{ we have
}\card{A_\delta\setminus\name{\tau}^\ast}=\aleph_0".
\]
We fix
in $V$ a continuously increasing sequence $\bar{N}=\langle
N_i:\,i<\omega_1\rangle$ of countable elementary submodels
of $\langle \HH(\chi),\in,<^\ast_\chi,p,\name{\tau},\bar{Q}\rangle$
such that $p^{N_0}=p^\ast$,while $\name{\tau}^{N_0}
=\name{\tau}^\ast$ and $\bar{Q}^{N_0}$ is our iteration $\bar{Q}$.
In addition, $\bar{N}\rest(i+1)\in N_{i+1}$ for all $i$.
For every $\gamma<\omega_1$, we can apply Claim \ref{2.5} to
$\bar{N}\rest(\gamma+1)$. Using this, we
can easily conclude that the set
\[
\begin{array}{lll}
C\deq\{\delta<\omega_1:\,  & (a) \,\, N_\delta\cap\omega_1=\delta\\
&(b)\,\,\,\,\delta\mbox{ is a limit ordinal}\\
&(c) \,\,\,\,\bar{\beta}_{\bar{N}\rest\delta,p^\ast,\name{\tau}^\ast}
\mbox{ and } \bar{r}_{\bar{N}\rest\delta,p^\ast,\name{\tau}^\ast}
\mbox{ are defined}\}
\end{array}
\]
is a club of $\omega_1$.
Let $\delta\in C$ be such that sequences $\bar{N}\rest\delta$ and
$\langle N_i^\delta:\,i<\delta\rangle$ have the same isomorphism type.
Let this be exemplified by $F=
\langle f_i:\,i<\delta\rangle$,
an increasing sequence of isomorphisms $f_i:\,N_i\into N^\delta_i$.
By our choice of constant symbols,
we also have that $F(\bar{Q})=\bar{Q}$, $F(p^\ast)=p^{N^\delta_0}$ and
$F(\name{\tau}^\ast)=\name{\tau}^{N^\delta_0}$.
By Subclaim \ref{X}, we have that
$\bar{\beta}_{\bar{N}^\delta,p^{N^\delta_0},\name{\tau}^{N^\delta_0}
}=
\bar{\beta}_{\bar{N}\rest\delta,p^\ast,\name{\tau}^\ast}$,
and
$\bar{r}_{\bar{N}^\delta,p^{N^\delta_0},\name{\tau}^{N^\delta_0}}
=F(\bar{r}_{\bar{N}
\rest\delta,p^\ast,\name{\tau}^\ast})$. We now let $\langle
\beta_n:\,n<\omega\rangle\deq
\bar{\beta}_{\bar{N}^\delta,p^{N^\delta_0},\name{\tau}^{N^\delta_0}}$.
By the definition of $\bar{r}$ and $\bar{\beta}$, there is $n_0$ and
condition
$q$ such that $q\forces`` \beta_n\in\name{\tau}^\ast"$ for all $n\ge
n_0$, and $q\ge p^\ast$. Hence
$q\forces``\card{A_\delta\setminus\name{\tau}^\ast}<\aleph_0"$,
which is in contradiction with the fact that $q\ge p^\ast$.

$\eop_{\ref{conc}}$
\end{Proof of the Claim}
$\eop_{\ref{1 not 0}}$

\end{Proof}

\begin{Note} 
(1) We note that the present result clearly implies that
$\clubsuit$ and $\stick$ are not the same (even without $CH$).

Clearly, $V^{P_{\omega_2}}\satisfies 2^{\aleph_0}=\aleph_2$. One of
the
ways to see this is to notice that under $CH$ the full $\clubsuit$ and
$\clubsuit^1$ agree (while $V^{P_{\omega_2}}\satisfies
2^{\aleph_0}\le\aleph_2$
obviously).

(2) Note that the sequence $\langle A_\delta:\,\delta<\omega_1\rangle$
exemplifying $\clubsuit^1$ in $V^P$, is in fact a sequence in $V$. 
\end{Note}

For clarity of
presentations we decided to give details of the proof of Theorem \ref{1 not 0}
rather than Theorem \ref{1 not bullet} below,
which is of course stronger than Theorem \ref{1 not 0}.
Now the obvious changes to the proof of Theorem \ref{1 not 0}
(just change the definition of $\name{Q}_\beta$) give

\begin{Theorem}\label{1 not bullet}  $CON(\clubsuit^1+\neg \clubsuit^\bullet)$.
\end{Theorem}

In the next section we encounter another similar proof, where the
changes needed to the proof of Theorem \ref{1 not 0}
are more significant, and we spell them out.

\section{Consistency of $\clubsuit^\bullet$ and $\neg\clubsuit^1$}\label{three}
\begin{Theorem}\label{bullet not 1} $CON(\clubsuit^\bullet+\neg\clubsuit^1)$.
\end{Theorem}

\begin{Proof} The proof is a modification of the proof from \S\ref{vec},
so we shall simply explain the changes, keeping all the non-mentioned
conventions and definitions in place.

Our iteration is again called $\bar{Q}=\langle P_\alpha,\name{Q}_\beta:\,
\alpha\le\omega_2,\beta<\omega_2\rangle$, but $\name{Q}_\beta$ will be redefined
below.

\begin{Definition}\label{ordern}
(1) A candidate for a $\clubsuit^1$ is a synonym for a candidate for $\clubsuit$.

(2) Suppose that $\beta<\omega_2$, and let us define
$Q_\beta$, while
working in $V^{P_\beta}$. It is defined the same way as in Definition 
\ref{order}(3), but we change the condition 1.$(iv)$ into
\[
(iv')\Dom(f)\cap A^\beta_\delta\mbox{ infinite }\implies
(\exists^\infty\gamma\in\Dom(f)\cap A^\beta_\delta)[f(\gamma)=0].
\]
\end{Definition}

\begin{Note} The following still hold with the new definition of the iteration
\begin{description}
\item{(1)} Claim \ref{ccc}.
\item{(2)} Conclusion \ref{concll}.
\item{(3)} Claim \ref{induction}.
\end{description}
\end{Note}

[Why? The same proofs.]

\begin{Claim}\label{bookn} It is possible to arrange the bookkeeping, so
that
$\forces_{P_{\omega_2}}\neg\clubsuit^1.$
\end{Claim}

\begin{Proof of the Claim}
It suffices to prove that for
every $\alpha^\ast<\omega_2$, in $V^{P_{\alpha^\ast}}$
we have
\[
\forces_{Q_{\alpha^\ast}}``\langle
A^{\alpha^\ast}_\delta:\,\delta<\omega_1\rangle
\mbox{ is not a }\clubsuit^1\mbox{-sequence}."
\]
Let $G$ be $Q_{\alpha^\ast}$-generic
over $V^{P_{\alpha^\ast}}$, and let $F\deq\bigcup G$. Let
$A\deq F^{-1}(\{1\})$. Suppose that $\card{A^{\alpha^\ast}_\delta\setminus
A}<\aleph_0$. We can find $p^\ast\in G$ which forces this, in fact without
loss of generality for some $\varepsilon<\delta$ we have
\[
p^\ast\forces``A^{\alpha^\ast}_\delta\setminus \name{A}\subseteq \varepsilon".
\]
But consider
\[
\II\deq\{q\ge p^\ast:\,(\exists \gamma\in (A^{\alpha^\ast}_\delta\setminus
\varepsilon)\cap\Dom(q))[q(\gamma)=0]\}.
\]
This set is dense above $p^\ast$: if $r\ge p^\ast$ is such that $\Dom(r)\cap
A^{\alpha^\ast}_\delta$ is infinite, then $r\in \II$. Otherwise,
let $\gamma=\min(A^{\alpha^\ast}_\delta\setminus (\Dom(r)\cup\varepsilon))$
and let $q\deq r\cup\{(\gamma, 0)\}.$ Contradiction.
$\eop_{\ref{bookn}}$
\end{Proof of the Claim}

\begin{Definition}\label{badn}
 Suppose that
\begin{description}
\item{(a)} $\gamma<\omega_1$,
\item{(b)} $\bar{N}=\langle N_i:\,i\le\gamma\rangle$ is a continuous
increasing
sequence of countable elementary submodels of $\langle
\HH(\chi),\in,<^\ast_\chi\rangle$,
\item{(c)} $\name{\tau},\bar{Q}\in N_0$ and $p\in P_{\omega_2}\cap
N_0$, 
\item{(d)} $p\forces``\name{\tau}\in [\omega_1]^{\aleph_1}"$ and
\item{(e)} $\bar{N}\rest(i+1)\in N_{i+1}$ for $i<\gamma$.
\end{description}

We say that $\varepsilon\le\gamma$ is {\em bad}
for $(\bar{N},\name{\tau},p,\bar{Q})$ if $\varepsilon$ is a limit ordinal,
and there is no $m(\varepsilon)=
m({\bar{N}\rest\varepsilon,p,\name{\tau}})<\omega$ and sequences
$\langle r^m_n:\,n<\omega\rangle$ and $\langle \beta^m_n:\,n<\omega\rangle$
for $m\le m(\varepsilon)$ such that $r^m_n, \beta^m_n\in N_\varepsilon$
and
\begin{description}
\item{(1)} $ r^m_n\forces_{P_{\omega_2}}``\beta_n^m\in \name{\tau}$",
\item{(2)} $\bigcup_{n\in \omega}\beta_n^m=N_\varepsilon\cap\omega_1,$
\item{(3)} $ r_n^m\ge p \mbox{ for all } n,m,$
\item{(4)} $\beta_n^m \mbox{ increase with }n,$ 
\item{(5)} $\mbox{ for some }m\le m(\varepsilon)$
the set $\{r_n^m:\,n<\omega\}$
has an upper bound in $P_{\omega_2}$
\item{(6)} $m(\varepsilon)$ and
${\bar{r}}_{\bar{N}\rest\varepsilon,p,\name{\tau}}
\deq\left\langle\langle r_n^m:\,n<\omega\rangle:\,m<m(\varepsilon)
\right\rangle$ and 
${\bar{\beta}}_{\bar{N}
\rest\varepsilon,p,\name{\tau}}\deq
\left\langle\langle \beta_n^m:\,n<\omega\rangle:\,m<m(\varepsilon)
\right\rangle$
are definable in $(\HH(\chi)^V,\in,<^\ast_\chi)$ from the isomorphism
type of
$(\langle N_\xi:\,\xi\le\varepsilon\rangle,
p,\name{\tau},\bar{Q})$ (we shall sometimes
abbreviate this by saying that these objects are defined in a
{\em canonical} way).
\end{description}
\end{Definition}

\begin{Main Claim}\label{2.5n} Suppose that $\bar{N},\gamma,p$ and
$\name{\tau}$
are as in Definition \ref{badn}.
\underline{Then} the set
\[
B\deq\{\varepsilon\le\gamma:\,\varepsilon\mbox{ bad for
}(\bar{N},\name{\tau},p,\bar{Q})\}
\]
has order type $<\omega^\omega$.
\end{Main Claim}

\begin{Proof of the Main Claim} Fix such $\bar{N},\gamma,p$ and
$\name{\tau}$.
We define $\bar{p}=
\bar{p}(\gamma,\bar{N},\name{\tau}, p)$ as in Subclaim \ref{p} and
$p_\ast, q^\ast, w^\ast, b$ as in the proof of Main Claim \ref{2.5}.
We shall show that $B\subseteq b$, by taking
any limit ordinal $\varepsilon^\ast\in (\gamma+1)\setminus b$ and showing
that it is not in $B$. 

Given $\varepsilon^\ast$,
we define $n^\ast$, $\xi$
and $\langle r_n:\,n<\omega\rangle$ and $\langle\beta_n:\,
n<\omega\rangle$the way we did in the proof of Main Claim \ref{2.5}.
We let $m(\varepsilon^\ast)=2^{n^\ast}-1$. For $m\le m(\varepsilon^\ast)$, we let
$\{i^m_n:\,n<\omega\}$ be the increasing enumeration of
$\{i<\omega:\,i=m(\mbox{mod }2^{n^\ast})\}$ and let $r^m_n=r_{i^m_n}$
and $\beta^m_n=\beta_{i^m_n}$. 
We shall show that
for some $m\le m(\varepsilon^\ast)$, the
sequence $\langle r_n^m:\,n<\omega\rangle$ has an upper bound in $P_{\omega_2}$.
Recall the definition of $\alpha_l$ for $l\le n^\ast$ from the proof of Main
Claim \ref{2.5}.
Notice that it is not a priori clear that $\bigcup_{n<\omega}r_n^m$
is a condition, as it may happen that for some $l<n^\ast$
it is forced that
$\name{X}_l\deq\bigcup_{n<\omega}\Dom(r^m_n)(\alpha_l)\cap
{\name{A}}^{\alpha_l}_{N_{\varepsilon^\ast}\cap\omega_1}$ is infinite, yet
$\bigcup_{n<\omega}r^m_n(\alpha_l)\rest\name{X}_l$ is 0 only finitely
often.

By induction on
$l\le n^{\ast}$
we choose $q^\ast_l\in P_{\alpha_l}$ and $k_l<2^l$, so that
\begin{description}
\item{$(a)$} $q_l^\ast\ge p_\ast\rest\alpha_l$,
\item{$(b)$} $(\forall n<\omega)[n=k_l(\mbox{mod }2^l) \implies
r_n\rest\alpha_l\le q_l^\ast]$.
\item{(c)} $q^\ast_l\le q^\ast_{l+1}$.
\end{description}
This clearly suffices, as we have that $q_{n^\ast}\in P_{\omega_2}$
is a common upper bound of $\{r^{k_{n^\ast}}_n:\,n<\omega\}$.

Let $q^\ast_0\deq q^{\ast}\rest
\alpha_0=
p_\ast\rest \alpha_0$.

Given $q_l^\ast\in P_{\alpha_l}$ and $k_l<2^l$ for some
$l<n^{\ast}$. Let $\Gamma\deq\{n<\omega:\,n=k_l(\mbox{mod }2^l)\}$.
Let $k'_1\deq k_l$ and $k'_2\deq k_l+2^l$.
Then $\Gamma=\Gamma_1\cup\Gamma_2$, where $\Gamma_1$ and $\Gamma_2$ are
infinite disjoint and defined by the following, for $j\in \{1,2\}$.
\[
\Gamma_j\deq\{n\in \Gamma:\,n=k'_j(\mbox{mod }2^{l+1})\}.
\]
If 
\[
q^\ast_l\forces``\bigcup_{n\in\Gamma_j}\Dom(r_n(\alpha_l))\bigcap
\name{A}^{\alpha_l}_{N_{\varepsilon^\ast}\cap\omega_1}\mbox{ finite}"
\]
for at least one $j\in\{1,2\}$,
let $j^\ast$ be the smallest such $j$ and let $k_{l+1}\deq k'_{j^\ast}$.
Let
\[
q^\ast_{l+1}\deq q^\ast_l\frown\{(\alpha_l,\bigcup_{n\in \Gamma_{j^\ast}}
r_n(\alpha_l))\}\frown p_\ast\rest(\alpha_l,\alpha_{l+1}).
\]
Otherwise,
we can find some $q'_l\in P_{\alpha_l}$ such that $q'_l\ge q^\ast_l$ and
\[
q'_l\forces ``\bigcup_{n\in\Gamma_2}\Dom(r_n(\alpha_l))\bigcap
\name{A}^{\alpha_l}_{N_{\varepsilon^\ast}\cap\omega_1}\mbox{ infinite}".
\]
Let $j^\ast\deq 1$ and $k_{l+1}\deq k'_1$, and let
\[
q^\ast_{l+1}\deq q'_l\frown\{(\alpha_l,\bigcup_{n\in \Gamma_{1}}
r_n(\alpha_l)\bigcup 0_{\bigcup_{n\in\Gamma_2}\Dom(r_n(\alpha_l))\setminus\xi}
)\}\frown p_\ast\rest(\alpha_l,\alpha_{l+1}).
\]
(Remember that for $n_1\neq n_2$, we have that $\Dom(r_{n_1}(\alpha_l))
\setminus\xi$ and $\Dom(r_{n_2}(\alpha_l))
\setminus\xi$ are disjoint.)

Observe, similarly to Subclaim \ref{X},
that the choice of $\bar{r}$ and $\bar{\beta}$ in this proof was canonical.
$\eop_{\ref{2.5n}}$
\end{Proof of the Main Claim}

\begin{Claim}\label{concn} $\forces_{P_{\omega_2}}\clubsuit^\bullet$.
\end{Claim}

\begin{Proof of the Claim}
Let
$\left\langle \bar{N}^\delta=\langle N^\delta_i:\,i<\delta\rangle:
\,\delta<\omega_1\right\rangle$ be as in the proof of Claim \ref{conc},
as well as $N^\delta$ for
limit ordinal $\delta<\omega_1$.

For limit $\delta<\omega_1$, we define $n^\ast(\delta)$ and
$\langle A_\delta^m:\,m\le m^\ast(\delta)\rangle$ as follows.

If
$\bar{\beta}_{\bar{N}^\delta,p^{N_\delta},\name{\tau}^{N_\delta}}$
and $\bar{r}_{\bar{N}^\delta,p^{N_\delta},\name{\tau}^{N_\delta}}$
are well
defined, then we let $m^\ast(\delta)\deq
m_{\bar{N}^\delta,p^{N_\delta},\name{\tau}^{N_\delta}}$ and for
$m\le m^\ast(\delta)$ we let $A_\delta^m\deq\{\beta^m_n:\,n<\omega\}$.
Otherwise, we let $m^\ast_\delta=0$ and
$A_\delta^0$ be the range of any cofinal $\omega$-sequence in $\delta$.

We claim that
\[
\left\langle \langle A_\delta^m:\,m\le m^\ast(\delta)\rangle
:\,\delta<\omega_1\right\rangle
\]
exemplifies that $V^P\models\clubsuit^\bullet(\omega_1)$.

Suppose that
\[
p^\ast\forces``\name{\tau}^\ast\in [\omega_1]^{\aleph_1}\mbox{ and for
all
}\delta, m\mbox{ we have
}A_\delta^m\setminus\name{\tau}^\ast\neq\emptyset".
\]
Let $\bar{N}$, $C$, $\delta$ and $F$ be as in the proof of Claim \ref{conc}.
It is easily seen that $q_{n^\ast}$ obtained as in the proof
of Main Claim \ref{2.5n} exemplifies a contradiction.
$\eop_{\ref{concn}}$
\end{Proof of the Claim}
$\eop_{\ref{bullet not 1}}$
\end{Proof}

\end{document}